\documentclass{notices}
\usepackage{amsfonts,amssymb,amsmath,amscd}

\usepackage{graphicx}
\usepackage{xcolor}
\usepackage{wrapfig}
\usepackage{CJKutf8}
\usepackage[T1]{fontenc}
\usepackage{hyperref}

\definecolor{Bittersweet}{rgb}{1.0, 0.44, 0.37}
\definecolor{KateColour}{rgb}{0.3, 0.6, 0.3}%56, 156, 64

\newcommand{\KS}[1]{}{}

\newcommand{\EH}[1]{}{}

\newcommand{\RR}{\mathbb{R}}

\title{Mathematical Illustration: Proof, Intuition and AI}
\author{
   Edmund Harriss
   \affil{
        Department of Mathematical Sciences and School of Art, University of Arkansas
   }
   \and
   Katherine E.~Stange
   \affil{
        Department of Mathematics, University of Colorado Boulder
   }
}

\date{\today}

\begin{document}

\maketitle

\section{Understanding}

``Drawing was knowing,'' explained Fr\'ed\'eric Brechenmacher, speaking in the Amphitheatre Charles Hermite.  He traced the history of mathematical education in France, beginning with the 18th century establishment of French royal engineering schools, through to the rising independence of pure mathematics as its own subject in the turn of the 20th century \cite{brechenmacher2026flat}. 

Throughout this time in France, the subject of geometry at the highest levels was taught in its largest proportion through \emph{model drawing}, not lecture.  Henri Poincar\'e documented his struggles to pass the model drawing requirements for the entrance to the \'Ecole Polytechnique.  Gaspard Monge, who designed the original curriculum at the school, wrote 

\begin{quote}
    the student who will have trained simultaneously his intelligence and his skills with his hands, will get, as the price of his double work, the exact description of the knowledge he will have acquired.~\cite{brechenmacher2026flat}
\end{quote}

The authors were in the audience as part of a workshop on \emph{Rigorous Illustration in Mathematics}, held at the Institute Henri Poincar\'e in early 2026.  We were in attendance as co-organizers of a larger trimester on \emph{Illustration as a Mathematical Research Technique}.  We had previously been part of a group who argued for the importance of mathematical illustration in the Notices~\cite{Coulon:IIMR}. We were still learning. Brechenmacher was speaking as a historian of mathematics who works with the model collection at Insitute Henri Poincar\'e.  His talk was a reminder that there has not always been the same answer to the nature of mathematical knowledge.

\begin{figure}
    \centering
    \includegraphics[width=\linewidth]{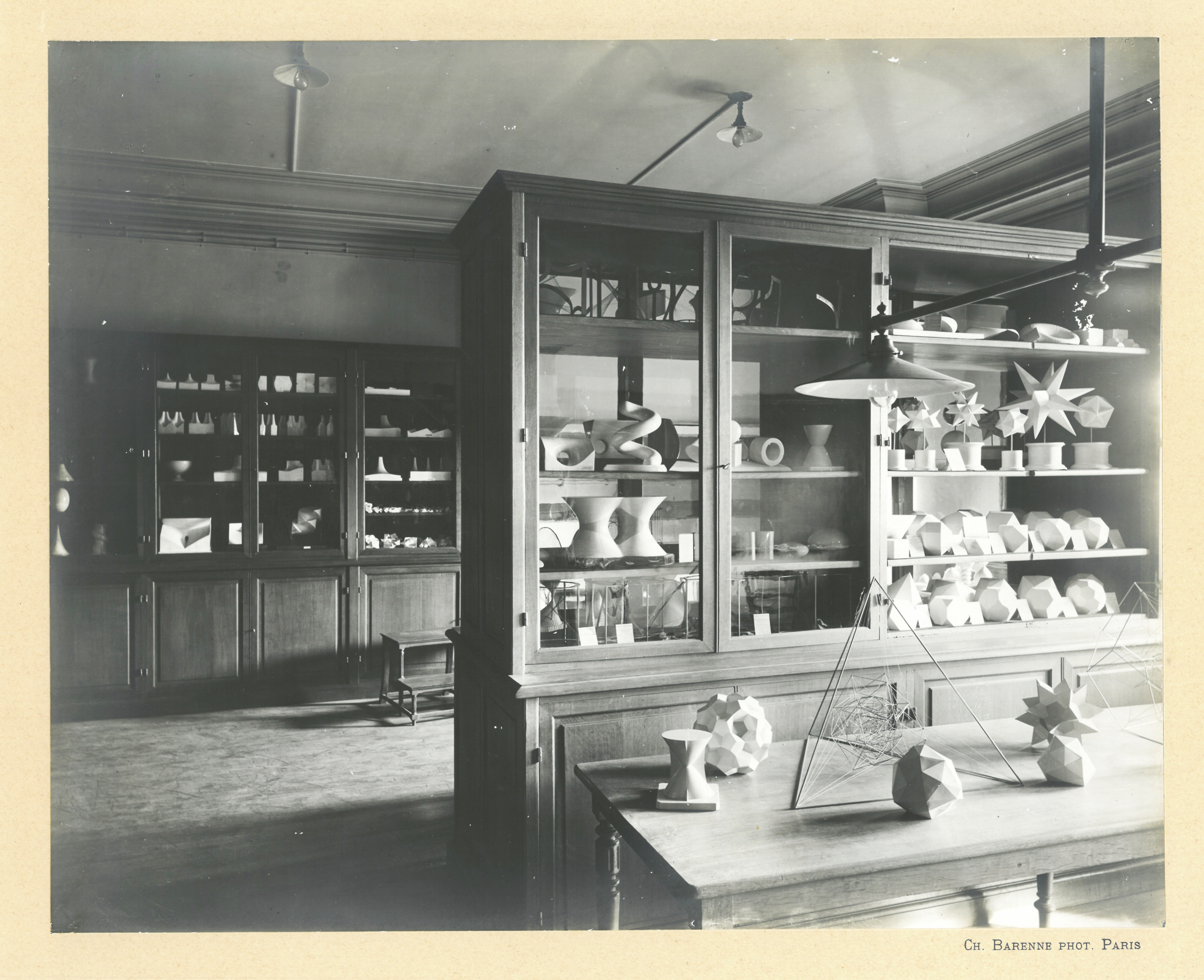}
    \caption{The models cabinet at the Laboratory of Higher Geometry (Faculty of Science, University of Paris), Charles Barenne, photographer, Paris, ca. 1910}
    \label{PoincareModels}
\end{figure}

The illustration of mathematics, by which we mean its incarnation in forms beyond text, be it drawing, graphics, sculpture or sound, among others\footnote{Vision being a dominant sense, throughout the article our discussion often focuses on images.  The same considerations apply, of course, to any of a much wider range of forms.},  has long held a position as a path to knowledge.  In Monge's curriculum, it was essential not only to view illustrations of geometry, but also to create them, to train one's hand and eye, and therefore learn something essential, something that was not encoded textually anywhere\footnote{
Monge's curriculum is a distinct historical thread to the development of models for cutting edge research. The German tradition of Brill and Schilling is an example.}.  It was only in the early 1900's that mathematics turned course fully to embrace rigorous proof as the primary, almost exclusive, standard for knowing.  In the 1910's, the Sorbonne still had a \emph{cabinet de math\'ematiques} of plaster models (Figure~\ref{PoincareModels}).  But from that time to this, these models have slowly drifted into the dusty back cabinets of history.  They can now be viewed in the IHP library collection.

Now, we stand in the moment of the Leiden Declaration \cite{leidendeclaration2026}, forced to ask again, what is mathematical knowledge?  If an LLM provides a rigorously verifiable but incomprehensible or unenlightening proof of a famous conjecture, what have we gained?  

Among the values upon which the Leiden Declaration describes itself as resting, we read:

\begin{quote}
Mathematics produces not only a body of results, but also understanding, clarity, and judgment among the communities of mathematicians who have shaped them, often in the context of their own autonomously guided research. This expert knowledge is essential, both to effectively use mathematics, and to continue to articulate new and significant research questions. A key source of strength of the discipline has long been the autonomous shaping of the direction of research and the methods used to pursue it.~\cite{leidendeclaration2026}
\end{quote}

The Declaration goes on to say, ``These characteristics of mathematics as a subject matter are also compatible with understanding mathematics as a human practice, and its place in the world.''  The mathematical community, having made this Declaration\footnote{The Leiden Declaration has been signed by thousands of mathematicians, but debate about it continues.  It is not our goal here to argue \emph{for} or \emph{against} the declaration.  Instead, we use it as a symbol and primary source touching on the moment in time in which we write.}, now has to revisit a number of fundamental questions.   What is the role of the formalization of mathematics, both in its human form and its automated one?  What is the role and value of intuition?  What aspects of the human community undertaking of the understanding of mathematics are mechanizable?

These questions are far from new, having been raised before by, for example, changes in standards of proof, including computer proofs. Yet they are brought into sharper focus with current developments in AI. That discussion is (and should be) wide-ranging, but in this article we more narrowly explore how illustration can be an important part of the answer\footnote{We recognize that aspects of this article could be out of date in a matter of months!  But this is a conversation we must drive forward.}.

Mathematicians have always drawn upon all their methods of understanding, including visual and spatial, tactile, and even emotional intuitions:  all the cognitive faculties that could be brought to bear. The community consensus on where the mathematical understanding lay, however, has been a delicate dance.

Our case is that illustration is a powerful tool for the practice of research mathematics.  In this article we encounter Euclid, Hilbert, Coxeter, Thurston, and philosophers of mathematical practice~\cite{DeToffoliMancosu2026} as we trace its historical traditions.  We wrestle with the current thinking in mathematics and philosophy.  We find that illustration has been part of the answer to what mathematics is and can be.  It deserves a thorough examination for its ability to carry rigorous meaning, and, especially in light of the discussion around AI, has great potential to reveal the wonder, beauty, utility and power of mathematics.

\section{The tradition of illustration}

Around the start of the twentieth century, the role of the models at the Sorbonne and in Monge's drawing classes were in decline in their impact on mathematical culture.  The deliberate sidelining of intuition, and with it illustration, was not careless. It was part of a broad, positive move towards formalism. The nineteenth century had seen dramatic growth in the scope of ideas mathematics could explore, and had revealed many surprises, such as non-euclidean geometry and continuous but nowhere differentiable functions. Such examples challenged the definitions of mathematics people were working with.  While some, like Poincar\'e, were concerned about mathematics being reduced to a formal game, the approach proved itself, with mathematics flourishing in the twentieth century into a cascade not just of new results but whole new areas. 

David Hilbert's axiomatics, through Bertrand Russell's foundations, to Bourbaki's systematic reconstruction, gave mathematics firm ground to build on:  formal proofs from an agreed-upon central collection of definitions and statements.  We are so used to mathematics being built this way today that we forget the staggering intellectual achievement it represents:  bringing the diverse collection of objects studied by mathematicians before the twentieth century into one cohesive framework. Even potential setbacks, like the work of Kurt Gödel, revealed both the importance and the value of such a program.  Formal mathematics delivered genuine insight into itself and the world.

The program of formalization came with an explicit rejection of figures and images, especially in the case of Bourbaki.  In fact, the very beauty and storytelling power of an image is key to the formalist critique. A figure can direct attention, illuminating its claims, but masking its distortions. A Euclidean diagram shows at best how the result works in a single case. The Italian school of algebraic geometry is often cited as the canonical example of the beauty of the ideas seducing mathematicians away from correctness.  For example in 1993, Arthur Jaffe and Frank Quinn write:

\begin{quote}
More recently in this century the ``Italian school'' of algebraic geometry did not avoid major damage: it collapsed after a generation of brilliant speculation\ldots in 1946 the subject was still regarded with such suspicion that Weil felt he had to defend his interest in it;~\cite{jaffequinn1993}
\end{quote}

Although it lost its institutional and cultural support, the practice of making images and models as a path to mathematical understanding persisted. Inspired by Poincar\'e and his work, Jacques Hadamard published ``Psychology of Mathematical Invention'' in 1945, sharing his own thinking, and interviewed many of his contemporaries, surveying well beyond the famous names. He discusses extensively how the thinking of mathematicians often takes place in visual or even kinesthetic terms. He also distinguishes the process of discovering mathematical ideas from the formal words needed to prove it~\cite{hadamard1945}. He points out explicitly that many mathematicians do not think directly in words or symbols.

Donald Coxeter went even further, working from physical mirrors and models his whole career, detailing symmetric structures that were critical in studying group theory, a space of pure symbol~\cite{roberts2006}. Group theory, in the form of the discovery of a sporadic simple group, also launched the career of John Conway, who consciously made games, play, and models part of all his investigations. His office in Princeton became a temple to mathematical objects~\cite{roberts2015}. 

\begin{figure}
    \centering
    \includegraphics[width=\linewidth]{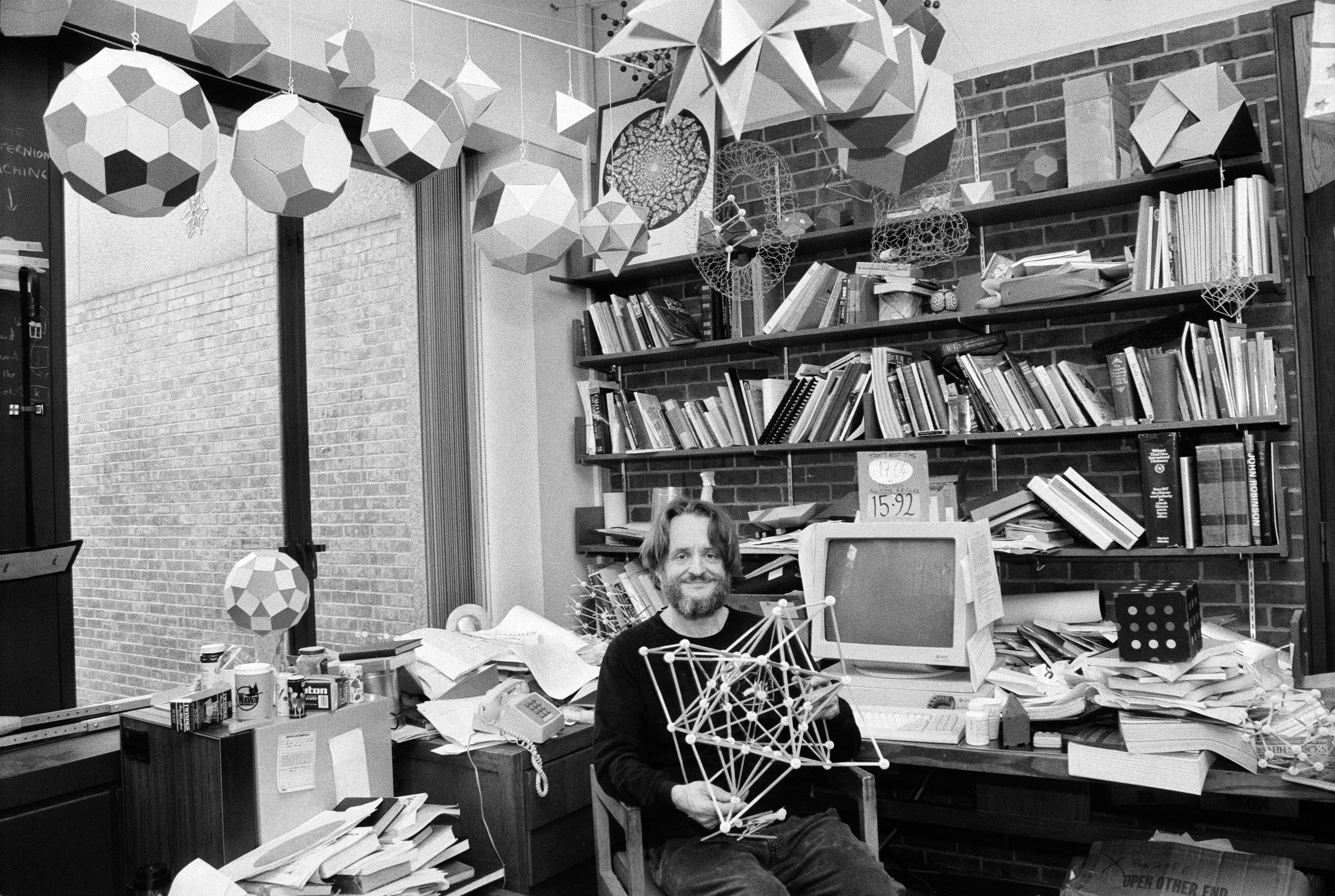}
    \caption{John Conway in his office at Princeton University in 1993. Photograph: Dith Pran/The New York Times}
\end{figure}

In the late 1980's Conway worked with Bill Thurston on the course ``Geometry and the Imagination,'' that took its title from the book by Stefan Cohn-Vossen and formalist hero Hilbert~\cite{Hilbert1952}. Conway and Thurston declared in the course notes that: \EH{Another removable citation}

\begin{quote}
Imagination, an essential part of mathematics, means not only the facility which is imaginative, but also the facility which calls to mind and manipulates mental images. One aim of the course is to develop the imagination.~\cite{conway2018}
\end{quote}

These beliefs were put into practice in the Geometry Center at the University of Minnesota, where the best documented version of the course was taught alongside Peter Doyle and Jane Gilman~\cite{conway2018}. The Geometry Center was active from 1991 to 1998, developing software and making videos. These videos, including ``Outside In'' and ``Not Knot'', have gathered millions of views on YouTube, even in today's competitive video market~\cite{outsidein,notknot,munznerarchive}. \KS{A place we could remove some references, as we are limited to twenty}

\begin{figure}
    \centering
    \includegraphics[width=\linewidth]{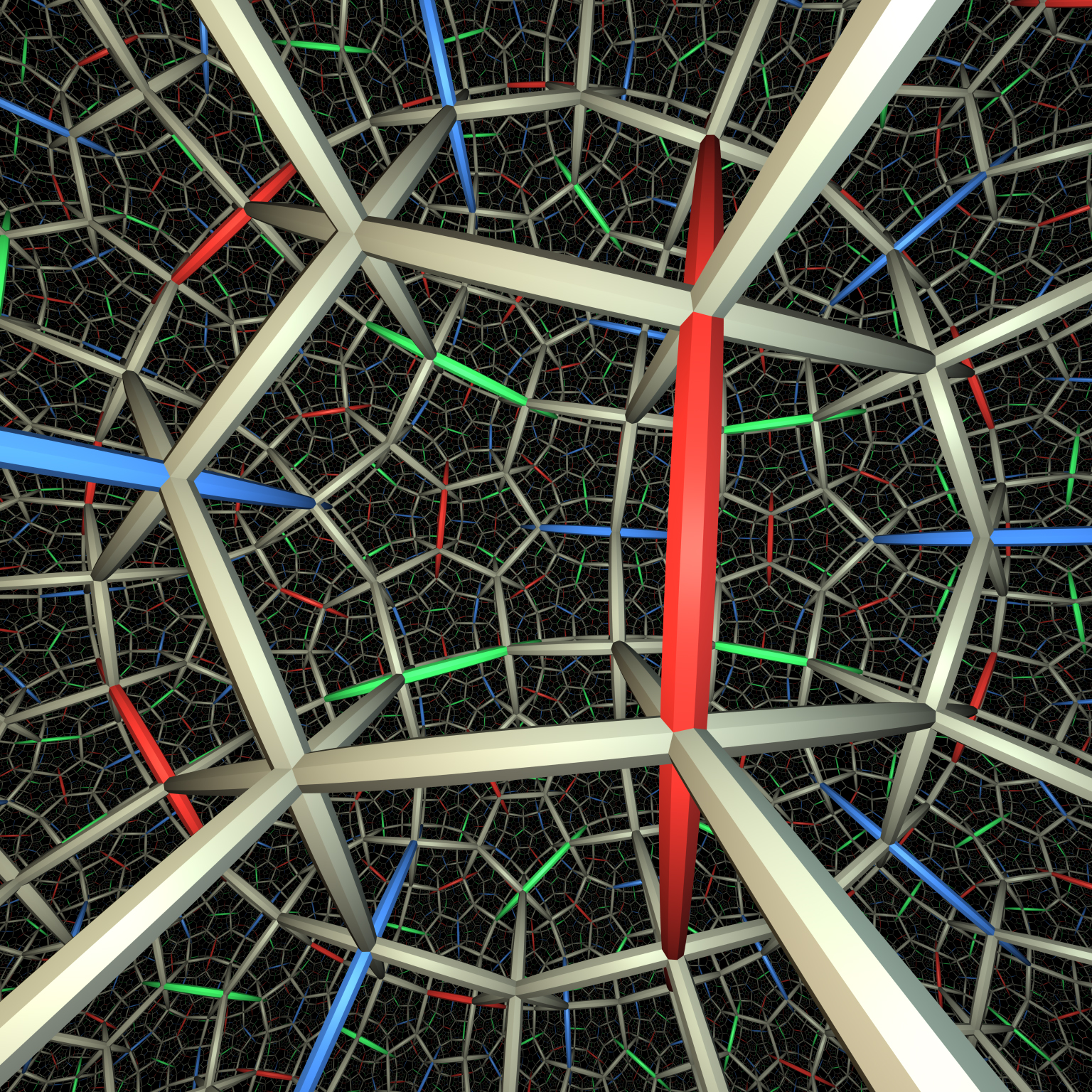}
    \caption{The poster for Not Knot video from the Geometry Center~\cite{notknot}, showing a view into hyperbolic space tiled by right-angle dodecahedra. © 1990 by The Geometry Center, University of Minnesota}
\end{figure}

Thurston's paper ``On Proof and Progress in Mathematics''~\cite{thurston1994} developed a powerful case for intuition and modeling in mathematical research, famously stating:

\begin{quote}
    Mathematicians usually have fewer and poorer figures in their papers and books than in their heads.~\cite{thurston1994}
\end{quote}

The paper also addresses what makes a theorem or a proof worthwhile to a mathematician:

\begin{quote}
    we should recognize that the humanly understandable and humanly checkable proofs that we actually do are what is most important to us, and that they are quite different from formal proofs.~\cite{thurston1994}
\end{quote}

Yet, this was not without context: it was written in response to the previously mentioned paper of Jaffe and Quinn~\cite{jaffequinn1993}, which had pushed against intuitive or speculative mathematics, and cited Thurston's geometrization as an example. They describe the result as:

\begin{quote}
    A grand insight delivered with beautiful but insufficient hints, the proof was never fully published. For many investigators this unredeemed claim became a roadblock rather than an inspiration.~\cite{jaffequinn1993}
\end{quote}

Thurston, therefore, was defending the integrity of his mathematics, a context which is often forgotten.  A followup paper had collected the responses of fifteen mathematicians, including Michael Atiyah, Armand Borel, Saunders Mac Lane and Karen Uhlenbeck~\cite{responses1994}.  This frequent  decontextualization is explicitly called out in Gila Hanna and Brendan Larvor's paper ``As Thurston says?'' \cite{hannalarvor2020}.  In fact, we have been following exactly the common pattern in mathematical writing Hanna and Larvor flag: telling a story through the famous names, in this case Hilbert and Poincar\'e, through Hadamard, to Coxeter and Conway, and Thurston. These are key figures, but they are voices in a wider community conversation that rarely reaches consensus.  Hanna and Larvor make this point by quoting Thurston themselves:

\begin{quote}
    Amongst other things, [Thurston's] article includes a version of the idea we pursue here, i.e., that we should not assume that, “there is uniform, objective and firmly established theory and practice of mathematical proof''.~\cite{hannalarvor2020}
\end{quote}

The history and philosophy of mathematical practice~\cite{DeToffoliMancosu2026} tries to find the commonalities and significant differences in those diverse individual accounts.  Here it demonstrates that illustration has been a consistent though often neglected part of the tapestry of mathematical practice in the twentieth century. Paolo Mancosu's 2008 collection, ``The Philosophy of Mathematical Practice''~\cite{mancosu2008} brought many of those working in the field together, finding:

\begin{quote}
    the epistemology of mathematics needs to be extended well beyond its present confines to address epistemological issues having to do with fruitfulness, evidence, visualization, diagrammatic reasoning, understanding, explanation, and other aspects of mathematical epistemology which are orthogonal to the problem of access to `abstract objects.'~\cite{mancosu2008}
\end{quote}

In that book, the case specifically for visual mathematics is made by Marcus Giaquinto. A more recent account by Giaquinto is given in the Stanford Encyclopedia of philosophy~\cite{giaquinto2020}. There, visual mathematics is described as `pervasive':  the text draws on examples from Euclidean geometry, knot theory, geometric group theory and others, addressing both strengths and challenges. It also concludes that

\begin{quote}
    By the late 20th century the mood had swung back in favor of visualization.~\cite{giaquinto2020}
\end{quote}

This story leaves out so much, but we will just mention a few more from the period:  David Mumford, Caroline Series and David Wright's research monograph on Kleinian groups, Indra's Pearls~\cite{Mumford2002}; computer geometry tools such as Cinderella~\cite{cinderella1999}; Tristam Needham's textbook on Visual Complex Analysis~\cite{Needham1997}; Daina Taimina's Crochet Hyperbolic planes giving physical form to hyperbolic geometry~\cite{Henderson2001};  Michelle Emmer's book, The Visual Mind, gave the fundamental connections between mathematics and art~\cite{Emmer1993}; and much more besides. An overview of the work during this period is sorely missing! \EH{I think these references stand for themselves and don't need the citations for the 20 limit}

\section{The rigor of illustration}

Despite that shift in mood, the Geometry Center itself did not survive. There has not been a full historical study of the causes, but a brief and near contemporary news item proposes a variety of explanations~\cite{mervis2002}, including that of David Dobkin, who chaired the center's board:
\begin{quote}
Its attempt to introduce computer graphics and visualization into pure mathematics and geometry was novel, and the community wasn't very welcoming. It wanted to change the field, but people weren't ready for that.~\cite{mervis2002}
\end{quote}

The mood might have swung in favor, but did not necessarily support a fundamental change, especially requiring a significant investment. The new tools of computer graphics created beautiful imagery and new ways of seeing mathematics without directly addressing the formalist rejection of illustration as part of intuition, and thus untrustworthy. The work of Hadamard and the claims of Conway and Thurston all push for the importance of imagination in building mathematical intuition. This is exactly the role to which formalism relegated illustration, tolerated as long as it was kept away from proof.  Russell's words from 1902, as he began work on the foundations of mathematics, were still echoing:

\begin{quote}
    It has been customary when Euclid, considered as a text-book, is attacked for his verbosity or his obscurity or his pedantry, to defend him on the ground that his logical excellence is transcendent,\ldots This claim, however, vanishes on a close inspection. \ldots A valid proof retains its demonstrative force when no figure is drawn, but very many of Euclid’s earlier proofs fail before this test.~\cite{russell1902}
\end{quote}

Russell did not have the last word on the subject, however. A key point in the development of the Philosophy of Mathematical Practice was Kenneth Manders' work~\cite{manders2008}, published in 2008, but having circulated as a manuscript since the 1990's. Manders carefully distinguished between the \emph{exact} and \emph{coexact} claims of a Euclidean proof. The \emph{coexact} claims are stable under small perturbation and so can be read off a diagram drawn with reasonable, but not exact precision. For example in the proof of Proposition 1 of Book 1 of Euclid, the intersection of the two circles to construct an equilateral triangle is coexact, so can be read off the diagram, but the equality of the three edges is exact, so must be proved syntactically, as it is in the Proposition.

\begin{figure}
    \centering
    \includegraphics[width=\linewidth]{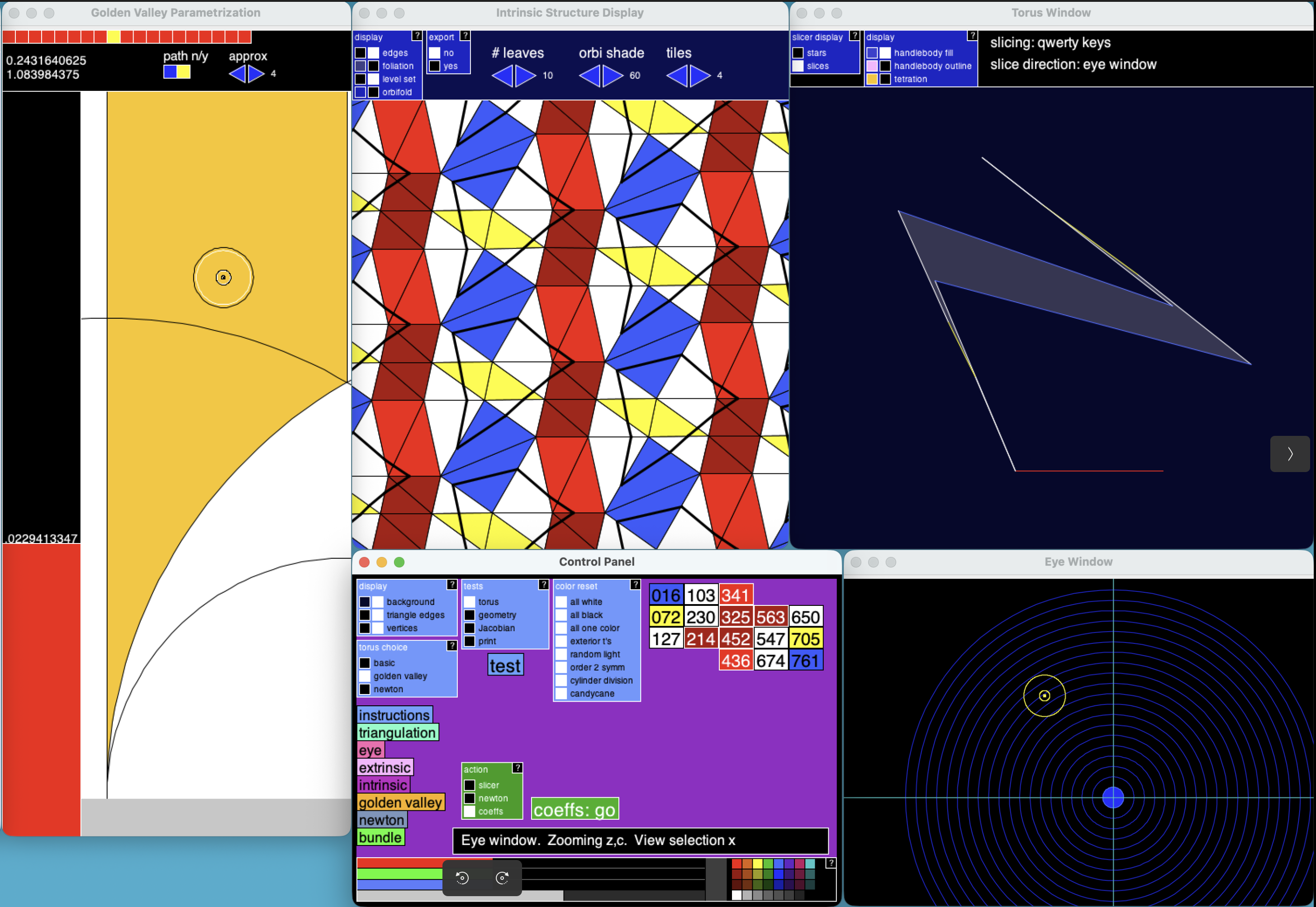}
    \caption{Screenshot from Rich Schwartz's research software, showing multiple illustrations of the same flat torus.}
\end{figure}

A precise version of the notion of coexact (where the perturbation is considered explicitly) is used in Rich Schwartz's discovery of the minimal flat torus. Using custom visualization software, he finds a model that is very close to the desired object by computer search and optimization. This model is robustly embedded, so it can be perturbed and stay embedded. Finally it is shown that within this radius of perturbation a model with the desired property will exist, by the implicit function theorem~\cite{schwartz2026}. 

Extending well beyond Euclid, the nature of diagrams has been well studied in recent decades, yet even in 2023, Sofia De Toffoli~\cite{detoffoli2023} feels the need to make a forceful case that they can carry weight in proof.  She begins,

\begin{quote}
Mathematical diagrams are frequently used in contemporary mathematics. They are, however, widely seen as not contributing to the justificatory force of proofs: they are considered to be either mere illustrations or shorthand for non-diagrammatic expressions.   Moreover, when they are used inferentially, they are seen as threatening the reliability of proofs.~\cite{detoffoli2023}
\end{quote}

Far from relegating illustration to the camp of intuition, unworthy of formalism, De Toffoli's exploration of diagram shows that image-based work can lie fully within formalism, having ``justificatory force.''  So diagram can play a formalist role.  However, other sorts of visual reasoning are still slightly dismissed here as ``mere illustration,'' playing their more familiar role supporting exploration and intuition building, as in Conway and Thurston.  

The main pushback against wider illustration as a source of understanding has always been that its reasoning can be misleading, yet the Leiden declaration points out that:

\begin{quote}
    Current automated techniques can produce plausible but unreliable (or even incorrect) arguments which are difficult to distinguish from correct mathematical proofs. This applies not only to informal arguments, but also to formalizations, where the difficulty lies in the translation between computer-encoded and human presentations of concepts.~\cite{leidendeclaration2026}
\end{quote}

This itself is not a new thing: human written proofs and theorem statements have long been seen to risk misleading the reader, or at least providing no insight beyond the result's correctness. The role of proof in enlightening---not merely convincing---the reader has been well-studied~\cite{dalessandro2019}.  The enlightening power of illustrations, and their rigorous meaning, deserve similar scrutiny.  Valeria Giardino's framing of notations and visualizations as tools for representational thought provides a way forward~\cite{giardino2018}. She analyses them as anchors for reasoning, whose value lies in their representational affordances: the manipulations, readings, and inferences a representation invites and supports. In her account, a well-designed tool is one that invites the user to explore in mathematically sound ways. This links back to the model drawing with which we began. Monge's students were not simply recording geometry: it was their interaction with the models that developed their mathematical understanding. Giardino and De Toffoli define the enhanced manipulative imagination, developed through these interactions with illustrations, paired with mathematical training~\cite{DeToffoliGiardino2014}.

The value of an image, like that of a proof, lies as much in the engagement it produces as in any property of the artifact alone. 

Fragment shader programming is an instructive, modern, example. A shader is a simple piece of code that takes as input a pixel coordinate and outputs a color to paint that pixel.  Shaders exploit the huge parallel computing power of the GPU, resulting in lightning-fast graphics. A shader program is a mathematical object in its own right, akin to an implicit equation for a graph, a function from the plane $\RR^2$ (of the screen) to a color space (generically a subset of $\RR^3$). To equate the results of implicit code with an explicit object requires a mathematical argument or proof. As an example, a simple shader computes the colour of a pixel in terms of its distance to the origin:  red if it is less than a unit.  This is the implicit code.  The observer sees a filled circle.  The code confirms the circle (up to the limits of implementation\footnote{Regarding this distinction between the code and its implementation.  A shader-coded circle running on a grid of rectangular pixels changes the norm of the space, resulting in an ellipse instead of a circle, albeit from the same code.  The relevance of implementation is not restricted to shaders: vector graphics software often uses splines to draw a circle, so that the result is not 45 degree rotationally symmetric (this can be checked in Adobe Illustrator 2025, by copying, rotating, and zooming in to see two distinct lines).  There is perhaps a sense in which the implementation layer creates its own set of axioms.}), as defined by the explicit equation $|z| \le 1$.

\begin{figure}
    \centering
    \includegraphics[width=\linewidth]{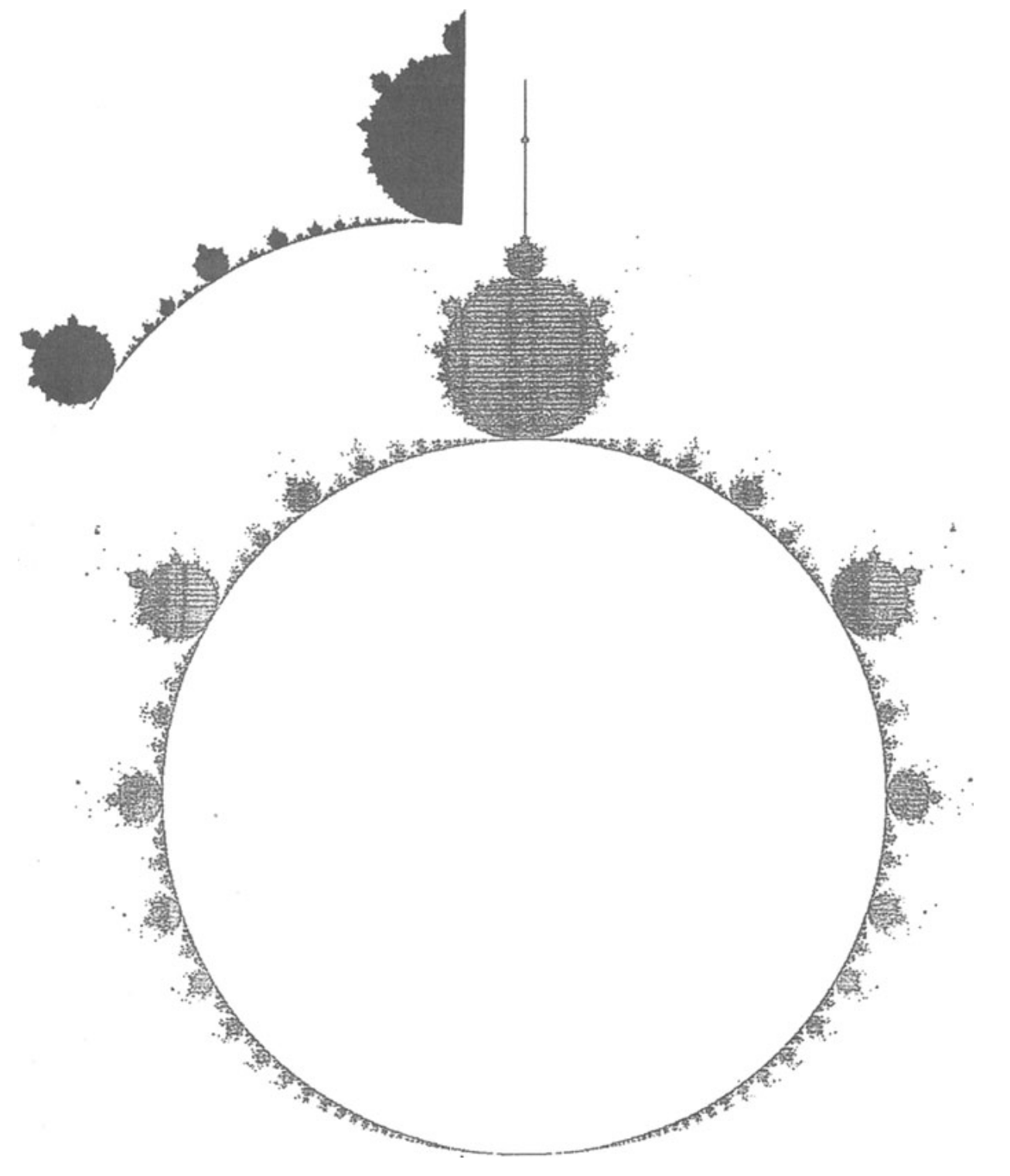}
    \caption{This reproduction of Figure C3-1 from Benoit Mandelbrot's book \emph{Fractals and Chaos} shows an early image of a Mandelbrot set (under a less standard parametrization).  The main image is Mandelbrot's original, motivating his observation in \cite{Mandelbrot1980Fractal} that it contained small specks that appeared disconnected.  In the upper left inset is a reproduction of how the image was `improved' by the publishers, much to the consternation of Mandelbrot, who later replaced the missing specks with hand-drawn ones when handing out reprints.}
    \label{fig:mandel}
\end{figure}

An only slightly more complicated algorithm gives the Mandelbrot set. Here the pixel colour is determined by whether an iterated quadratic function escapes a bounded region in a fixed number of steps. The result is an approximation of a mathematical object (since we have only a discrete spectrum of colors, a finite region of pixels, and the stopping test is an approximation of a limiting behavior). It is also a mathematical object in its own right, the product of the algorithm chosen. The resulting image can be studied on its own merits.  Here we might write theorems and proofs about the properties of the image in its role as an approximation.  In Giardino's terms, such a proof certifies which of the image's affordances can be trusted.  Early renderings of the Mandelbrot set revealed many unforeseen and insightful properties, but also afforded at least one false inference: the set appeared disconnected (Figure~\ref{fig:mandel}), but Adrien Douady and John Hubbard proved that it was connected~\cite{douadyhubbard1982}. 

\begin{figure}
    \centering
    \includegraphics[width=\linewidth]{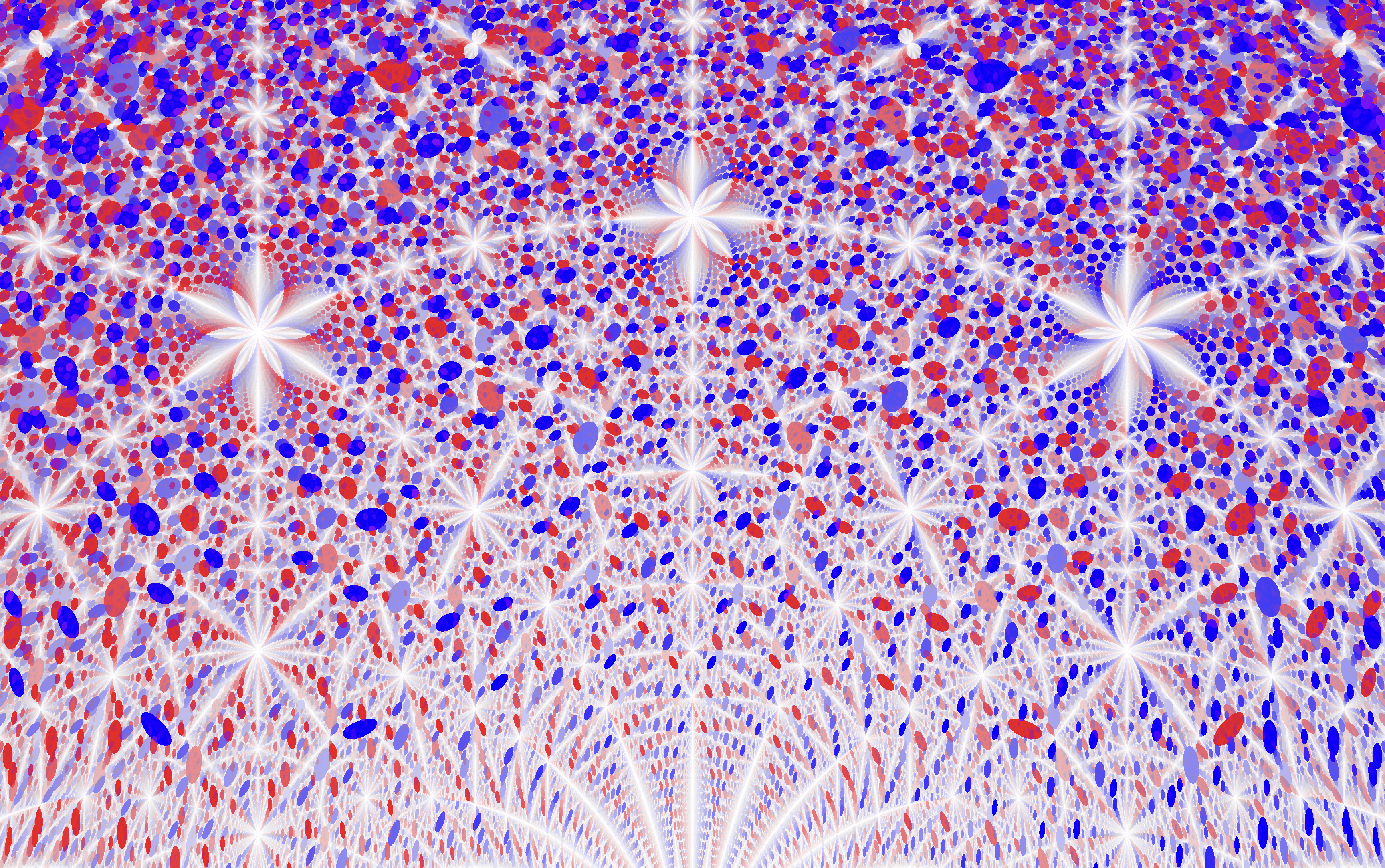}
    \caption{A shader rendering of cubic algebraic numbers in the complex upper half plane, developed by Steve Trettel as part of our joint work on algebraic starscapes.  The cubics are rendered as ellipses because the algorithm consists of a geometric transformation implemented via raytracing, from the lattice of spheres representing cubic polynomials in coefficient space. The code is available here: \url{https://edmund.mathematicians.org.uk/visualisers/cubic-shader/}}
    \label{ShaderCubics}
\end{figure}

Illustration has an uncanny ability to generate unexpected insights.  We, the authors, motivated by speed, have been reillustrating our past work on algebraic starscapes using shaders. Diophantine approximation became an essential technology for creating the images, rather than something we observed in them.  Inevitably---but nevertheless a pleasant surprise---this changed the fundamental mathematics.  We suddenly had a field of ellipses, not discs (Figure~\ref{ShaderCubics}).  Although the distortion of the discs was a casualty of the new algorithm, the new shapes add a new dimension of visual information about the underlying geometric transformations relating coefficients and roots.  The image itself became much more compelling, structural elements such as the relation to the unit circle became clearer, and in terms of its beauty.  Beauty is not a side-effect of illustration, any more than it a side-effect of proof.  Beautiful ideas and beautiful pictures are fellow travelers, and mathematicians can and do use both to guide intuition and exploration.  We are still trying to understand and make rigorous what we are seeing. 

A robust discussion of the practice of mathematical illustration must call on the work of other disciplines: mathematical education, which places meaning in a pedagogic setting; scientific illustration, which examines the faithfulness of illustration; and art, for notions of critique and the beholder's share. 

We began this article with a talk from the Rigorous Illustration workshop at IHP.  That workshop aimed to impel this interdisciplinary endeavor, and included historians, mathematicians and artists.  From mathematics education, Hortensia Soto's concluding talk offered a framework where engagement with physical materials supports learners' actions, perceptions, and sensations, leading to mathematical understanding~\cite{soto2026rigorous}. This paper had its origin in the discussions that came out of that workshop. 

\section{The potential of illustration}

We sit at one of the periodic moments of flux in the community, following the turns to formalism in the 1910s and to computing in the 1990s, when the definition of the discipline is called into question. Our formalist castle has provided both knowledge and understanding, always intertwined. Automated proof may decouple them: a result can now arrive certified but unilluminating. Yet most mathematicians, asked to choose between the knowledge and the understanding a proof provides, will not affirm one at the expense of the other. The goals are unchanged; the methods must adapt.

Returning to illustration provides one such adaptation, strengthening threads already woven deep into the history of mathematics, having kept a vibrant presence even in the face of active challenge from formalism. Formalism rejected the incompleteness of intuition, and associated illustration only with intuition. As we have seen, however, illustration plays many roles beyond an aid to intuition.  Diagrams convince, careful representations challenge that very intuition, visual objects are valid subjects of reasoning in their own right, and illustration carries knowledge not available in textual form.  The beauty and aesthetic choices in an illustration can inspire new directions of inquiry.  In terms of building knowledge, illustration therefore plays a role more akin to proof. 

Of course, no one aims to replace proof with illustration. Nor is the push for illustration new. Although the Geometry Center itself did not survive, its inspiration continued. Bill Goldman, a board member at the center, created the Experimental Geometry Lab at the University of Maryland in 2000, with Rich Schwartz. Anton Lukyanenko, one of the alumni of the EGL, together with Jayadev Athreya, founded Illinois Geometry Lab, creating a model that spread to a network of similar geometry and math, under the umbrella of Geometry Labs United.  Labs host ``vertically integrated groups of students and researchers advancing mathematics through experimentation, visualization, and public engagement.''~\cite{glu2018}. Schwartz and the ICERM scientific advisory board, encouraged a new group to form around Illustrating Mathematics. ICERM hosted a week-long event in 2016, followed by a semester program in 2019. The Illustrating Mathematics community went on to organize the trimester program at IHP~\cite{Coulon:IIMR} which hosted Brechenmacher's talk with which we opened this article.

In parallel, the Philosophy of Mathematical Practice has been building its historical accounts and conceptual tools with an intent to expand beyond a purely formalist epistemology. Similarly, research in mathematics education has explored the question of how understanding can be built for mathematics whose rigor is well established.

And today, the advent of AI can not help but drive further innovation.  In just the past year, the advent of powerful vibecoding has suddenly made computer visualization orders of magnitude faster, easier and more accessible.  AI tools hold immense promise for expanding the ways we interact with mathematics, and for revealing its beauty. 

The definition and values of mathematics come not from points of authority or clean philosophical argument, but from the conversations and decisions of individuals, small groups and large organizations. These communities did not grow up in response to large language models, or the Leiden declaration. They are new expressions of traditions that go back through mathematical history, and they are already providing one thread of the positive, powerful, mathematical, response the current moment calls for. They support the positive case for the richness of the mathematical experience, made forcibly by Rodrigo Ochigame, one of its authors, in a recent article in the New York Times~\cite{nyt2026leidenai}:

\begin{quote}
    Mathematics is a rich form of cultural expression with an ancient history, and I am not worried that any technology will ever render it obsolete. Its most precious aspects, such as the collective quest to understand beautifully intricate ideas, and to explore the limits of the human imagination, cannot ever be automated. What I am worried about is that a handful of corporations are mobilizing their vast financial resources to impose an impoverished view of mathematics\ldots
\end{quote}

Before concluding, we want to strengthen the Leiden declaration's calls found in ``Recommendations for individual mathematicians'' items $7$ ``Participate in public discourse'' and $9$ ``Welcome new contributors'':

\begin{quote}
    \dots we encourage mathematicians to seek opportunities to cooperate with and support other researchers and creative professionals \ldots 
    The growing intersection of artificial intelligence and mathematics continues to attract researchers from other disciplines. We welcome this broadening of our community and the range of skills and perspectives these contributors bring. 
\end{quote}

This point is critical. Mathematicians have a level of understanding of technical material, and can speak with disciplinary knowledge on the mathematics underlying large language models, diffusion based image generation and other technologies, cutting through aspects of the hype. This can be combined with understanding from sciences, the humanities and the arts to develop a picture of how these systems work, and the social and cultural issues they raise. This is not a moment that belongs only to our discipline and we have a powerful potential role to play in support of the academy as a whole. 

Illustration has a long history.  It has consistently emerged as a tool for both understanding and proof.  It is a fundamentally human aspect of mathematics.  It has a role to play in answering our current cultural upheaval.  If we have convinced you of these things, then what can you do?  

In practice, the definition of mathematics is constructed by the community itself:  what we value collectively, promote at conferences, accept in publications, and speak about with enthusiasm to one another.  The most lasting way to support a practice is through community and cultural recognition.  Do you hold illustration in general as valuable?  There are dozens of small ways that respect and recognition can be manifest.  Students can be encouraged when they are drawn to illustration as a method of understanding.  Authors can be encouraged to illustrate, and space can be made for publishing innovative illustration backed by proofs of its affordances.  New perspectives can be supported during funding and hiring.  We can hold in mind that drawing, now extended by many more creative tools, still is knowing.

\section*{Acknowledgments}

The authors are indebted to the Institute Henri Poincar\'e, which hosted both trimester and workshop, as well as the participants of both, who shared energy and ideas through many fruitful discussions.  We are grateful to Pierre Arnoux, Jayadev Athreya, Gabriel Dorfsman-Hopkins, Valeria Giardino, Rich Schwartz, Hortensia Soto, Andrew Sutherland and Glen Whitney for early feedback on a draft, and to Henri Duvillard, Siobhan Roberts and Rich Schwartz for help sourcing images.

The authors acknowledge the support of the Institut Henri Poincar\'e (UAR 839 CNRS-Sorbonne Universit\'e), and LabEx CARMIN (ANR-10-LABX-59-01).  The first author would like to thank the university of Aix-Marseille and Pierre Arnoux for their hospitality.  The second author would like to thank the Institut des Hautes Études Scientifiques (IHES) in Bures-sur-Yvette, France, for its hospitality.  The second author was supported by NSF DMS-2401580 and a Joan and Joseph Birman Fellowship from the American Mathematical Society.

\section*{AI Declaration}

The authors made use of large language models (specifically Claude) during the preparation of this manuscript for bibliographic assistance, structuring and proofreading. All text and arguments are the authors' own.

\bibliography{references}

@book{hadamard1945,
  author = {Hadamard, Jacques},
  title = {The Psychology of Invention in the Mathematical Field},
  publisher = {Princeton University Press},
  year = {1945}
}

@book{roberts2006,
  author = {Roberts, Siobhan},
  title = {King of Infinite Space: {D}onald {C}oxeter, the Man Who Saved Geometry},
  publisher = {Walker \& Company},
  year = {2006}
}

@book{roberts2015,
  author = {Roberts, Siobhan},
  title = {Genius at Play: The Curious Mind of {J}ohn {H}orton {C}onway},
  publisher = {Bloomsbury},
  year = {2015}
}

@article{thurston1994,
  author = {Thurston, William P.},
  title = {On Proof and Progress in Mathematics},
  journal = {Bulletin of the American Mathematical Society},
  volume = {30},
  pages = {161--177},
  year = {1994}
}

@article{jaffequinn1993,
  author = {Jaffe, Arthur and Quinn, Frank},
  title = {{`}{T}heoretical Mathematics': Toward a Cultural Synthesis of Mathematics and Theoretical Physics},
  journal = {Bulletin of the American Mathematical Society},
  volume = {29},
  pages = {1--13},
  year = {1993}
}

@article{responses1994,
  author = {Atiyah, Michael and Borel, Armand and Chern, Shiing-Shen and others},
  title = {Responses to {`{T}heoretical} {M}athematics: {T}oward a Cultural Synthesis of Mathematics and Theoretical Physics,' by {A}. {J}affe and {F}. {Q}uinn},
  journal = {Bulletin of the American Mathematical Society},
  volume = {30},
  pages = {178--207},
  year = {1994}
}

@article{mervis2002,
  author = {Mervis, Jeffrey},
  title = {The {G}eometry {C}enter, 1991--1998. {RIP}},
  journal = {Science},
  volume = {297},
  pages = {508},
  year = {2002}
}

@article{hannalarvor2020,
  author = {Hanna, Gila and Larvor, Brendan},
  title = {As {T}hurston says? {O}n using quotations from famous mathematicians to make points about philosophy and education},
  journal = {ZDM -- Mathematics Education},
  volume = {52},
  number = {6},
  pages = {1137--1147},
  year = {2020}
}

@book{mancosu2008,
  editor = {Mancosu, Paolo},
  title = {The Philosophy of Mathematical Practice},
  publisher = {Oxford University Press},
  year = {2008}
}

@incollection{manders2008,
  author = {Manders, Kenneth},
  title = {The {E}uclidean Diagram},
  booktitle = {The Philosophy of Mathematical Practice},
  editor = {Mancosu, Paolo},
  publisher = {Oxford University Press},
  pages = {80--133},
  year = {2008},
  note = {Originally circulated as manuscript, 1995}
}

@InCollection{giaquinto2020,
	author       =	{Giaquinto, Marcus},
	title        =	{{The Epistemology of Visual Thinking in Mathematics}},
	booktitle    =	{The {Stanford} Encyclopedia of Philosophy},
	editor       =	{Edward N. Zalta},
	howpublished =	{\url{https://plato.stanford.edu/archives/spr2020/entries/epistemology-visual-thinking/}},
	year         =	{2020},
	edition      =	{{S}pring 2020},
	publisher    =	{Metaphysics Research Lab, Stanford University}
}

@article{giardino2018,
  author = {Giardino, Valeria},
  title = {Tools for Thought: The Case of Mathematics},
  journal = {Endeavour},
  volume = {42},
  number = {2--3},
  pages = {172--179},
  year = {2018}
}

@article{detoffoli2023,
  author = {De Toffoli, Silvia},
  title = {Who's Afraid of Mathematical Diagrams?},
  journal = {Philosophers' Imprint},
  volume = {23},
  year = {2023}
}

@article{douadyhubbard1982,
  author = {Douady, Adrien and Hubbard, John H.},
  title = {Itération des polynômes quadratiques complexes},
  journal = {Comptes Rendus de l'Académie des Sciences, Série I Mathématique},
  volume = {294},
  pages = {123--126},
  year = {1982}
}

@article{dalessandro2019,
  author = {D'Alessandro, William},
  title = {Explanation in Mathematics: Proofs and Practice},
  journal = {Philosophy Compass},
  volume = {14},
  number = {11},
  year = {2019}
}

@article{schwartz2026,
  author = {Schwartz, Richard Evan},
  title = {The most efficient origami torus},
  journal = {Proceedings of the National Academy of Sciences},
  volume = {123},
  number = {21},
  pages = {e2523301123},
  year = {2026},
  month = {May},
  doi = {10.1073/pnas.2523301123}
}

@article{glu2018,
  author = {Athreya, Jayadev and Dumas, David and Goldman, William and Grigorian, Sergey and Guzman, Rosemary and Hieronymi, Philipp and Lawton, Sean and Lukyanenko, Anton and Tyson, Jeremy and Wilson, Aaron},
  title = {Geometry Labs United: An Invitation},
  journal = {Notices of the American Mathematical Society},
  volume = {65},
  number = {9},
  pages = {1088--1093},
  year = {2018}
}

@article{russell1902,
  author = {Russell, Bertrand},
  title = {The Teaching of {E}uclid},
  journal = {The Mathematical Gazette},
  volume = {2},
  number = {33},
  pages = {165--167},
  year = {1902}
}

@misc{notknot,
  author = {Gunn, Charlie and Maxwell, Delle},
  title = {Not Knot},
  type = {Film},
  publisher = {The Geometry Center / A K Peters},
  year = {1991}
}

@misc{outsidein,
  author = {Levy, Silvio and Munzner, Tamara and Phillips, Mark},
  title = {Outside In},
  type = {Film},
  publisher = {The Geometry Center / A K Peters},
  year = {1994}
}

@misc{munznerarchive,
  author = {Munzner, Tamara},
  title = {Geometry {C}enter Videos, Revisited},
  howpublished = {\url{https://www.cs.ubc.ca/~tmm/gc/}},
  note = {Digitized films, tapes, and supplementary materials}
}

@book{Mumford2002,
  author = {Mumford, David and Series, Caroline and Wright, David},
  title = {Indra's Pearls: The Vision of {F}elix {K}lein},
  publisher = {Cambridge University Press},
  year = {2002},
  isbn = {978-0-521-35253-6}
}

@article{Henderson2001,
  author = {Henderson, David W. and Taimina, Daina},
  title = {Crocheting the Hyperbolic Plane},
  journal = {The Mathematical Intelligencer},
  volume = {23},
  number = {2},
  pages = {17--28},
  year = {2001},
  doi = {10.1007/BF03026623}
}

@book{Needham1997,
  author = {Needham, Tristan},
  title = {Visual Complex Analysis},
  publisher = {Oxford University Press},
  year = {1997},
  isbn = {978-0198534471}
}

@book{Emmer1993,
  editor = {Emmer, Michele},
  title = {The Visual Mind: Art and Mathematics},
  publisher = {MIT Press},
  address = {Cambridge, MA},
  year = {1993},
  series = {Leonardo Book Series},
  isbn = {978-0-262-05048-7}
}

@book{Hilbert1952,
  author = {Hilbert, David and Cohn-Vossen, Stephan},
  title = {Geometry and the Imagination},
  publisher = {Chelsea Publishing Company},
  address = {New York},
  year = {1952},
  note = {Translated by P. Nem\'enyi. Originally published as \emph{Anschauliche Geometrie}, Springer, Berlin, 1932}
}

@misc{conway2018,
      title={Geometry and the Imagination in Minneapolis}, 
      author={John H. Conway and Peter G. Doyle and Jane Gilman and William P. Thurston},
      year={2018},
      eprint={1804.03055},
      archivePrefix={arXiv},
      primaryClass={math.DG},
      url={https://arxiv.org/abs/1804.03055}, 
}

@book{cinderella1999,
  author    = {Richter-Gebert, J{\"u}rgen and Kortenkamp, Ulrich H.},
  title     = {The Interactive Geometry Software Cinderella},
  publisher = {Springer},
  address   = {Berlin},
  year      = {1999}
}

@article{Coulon:IIMR,
	author = {Coulon, R\'{e}mi and Dorfsman-Hopkins, Gabriel and Harriss, Edmund and Skrodzki, Martin and Stange, Katherine E. and Whitney, Glen},
	doi = {10.1090/noti2839},
	fjournal = {Notices of the American Mathematical Society},
	journal = {Notices Amer. Math. Soc.},
	number = {1},
	pages = {105--115},
	title = {On the Importance of Illustration for Mathematical Research},
	url = {https://doi.org/10.1090/noti2839},
	volume = {71},
	year = {2024}}

@article{nyt2026leidenai,
  author        = {Roberts, Siobhan},
  journal       = {{The New York Times}},
  title        = {{AI} Is Challenging the Core Values of Mathematics, Researchers Warn in “{L}eiden {D}eclaration”},
  year         = {2026},
  month        = jun,
  day          = {2},
  url          = {https://www.nytimes.com/2026/06/02/science/ai-mathematics-leiden-declaration.html},
  note         = {Accessed: 2026-06-15}
}

@misc{brechenmacher2026flat,
  author       = {Brechenmacher, Fr{\'e}d{\'e}ric},
  title        = {The Flat and the Relief: Geometric Material Models in 19th Century {France}},
  year         = {2026},
  month        = jan,
  day          = {20},
  howpublished = {Audiovisual resource, Institut Henri Poincar{\'e}},
  note         = {Workshop: Rigorous Illustrations---Their Creation and Evaluation for Mathematical Research (2026 T1 WS1). \url{https://dx.doi.org/10.57987/IHP.2026.T1.WS1.007}},
  doi          = {10.57987/IHP.2026.T1.WS1.007},
  url          = {https://www.carmin.tv/en/collections/2026-t1-ws1-rigorous-illustrations-their-creation-and-evaluation-for-mathematical-research/video/the-flat-and-the-relief-geometric-material-models-in-19th-century-france},
}

@misc{soto2026rigorous,
  author       = {Soto, Hortensia},
  title        = {Rigorous and Flexible Illustrations: Inviting Realms of Possibilities},
  year         = {2026},
  month        = jan,
  day          = {22},
  howpublished = {Audiovisual resource, Institut Henri Poincar{\'e}},
  note         = {Workshop: Rigorous Illustrations---Their Creation and Evaluation for Mathematical Research (2026 T1 WS1). \url{https://dx.doi.org/10.57987/IHP.2026.T1.WS1.017}},
  doi          = {10.57987/IHP.2026.T1.WS1.017},
  url          = {https://www.carmin.tv/en/collections/2026-t1-ws1-rigorous-illustrations-their-creation-and-evaluation-for-mathematical-research/video/rigorous-and-flexible-illustrations-inviting-realms-of-possibilities},
}

@misc{leidendeclaration2026,
  author       = {Alper, Jarod and Barany, Michael and {Chavarri Villarello}, Alain
                  and Dahmen, Sander and Dean, Walter and Ganapathy, Karthik
                  and Harris, Michael and Holmes, David and Jamnik, Mateja
                  and Kelk, Steven and Kra, Bryna and Martin, Ursula
                  and Naskr\k{e}cki, Bartosz and Ochigame, Rodrigo
                  and Portegies, Jim and Schmitt, Johannes},
  title        = {Leiden Declaration on Artificial Intelligence and Mathematics},
  year         = {2026},
  month        = jun,
  day          = {2},
  howpublished = {Community declaration, endorsed by the International Mathematical Union},
  doi          = {10.5281/zenodo.20302944},
  url          = {https://leidendeclaration.ai/},
  note         = {\url{https://doi.org/10.5281/zenodo.20302944}},
}

@article{Mandelbrot1980Fractal,
  author  = {Mandelbrot, Benoit B.},
  title   = {Fractal Aspects of the Iteration of $z \mapsto \lambda z(1-z)$ for Complex $\lambda$ and $z$},
  journal = {Annals of the New York Academy of Sciences},
  volume  = {357},
  number  = {1},
  pages   = {249--259},
  year    = {1980},
  doi     = {10.1111/j.1749-6632.1980.tb29690.x}
}

@incollection{DeToffoliMancosu2026,
  author    = {De Toffoli, Silvia and Mancosu, Paolo},
  title     = {The Philosophy of Mathematical Practice},
  booktitle = {The Stanford Encyclopedia of Philosophy},
  editor    = {Zalta, Edward N. and Nodelman, Uri},
  edition   = {Summer 2026},
  publisher = {Metaphysics Research Lab, Stanford University},
  year      = {2026},
  howpublished = {\url{https://plato.stanford.edu/archives/sum2026/entries/mathematical-practice/}}
}

@article{DeToffoliGiardino2014,
  author    = {De Toffoli, Silvia and Giardino, Valeria},
  title     = {Forms and Roles of Diagrams in Knot Theory},
  journal   = {Erkenntnis},
  year      = {2014},
  volume    = {79},
  number    = {4},
  pages     = {829--842},
  doi       = {10.1007/s10670-013-9568-7}
}

\end{document}